\documentstyle[12pt,spacing]{article}
\begin{document}
\begin{center}
\today\\[10pt]
{\Large\bf Co-elementary Equivalence, Co-elementary Maps, and Generalized Arcs}
\\[20pt]
Paul Bankston\\
Department of Mathematics, Statistics and Computer Science\\ 
Marquette University\\
Milwaukee, WI 53233\\[20pt]
\end{center}
\begin{abstract}
By a {\bf generalized arc\/} we mean a continuum with
exactly two non-separating points; an {\bf arc} is a metrizable generalized arc.
It is well known that any two arcs are homeomorphic (to the real closed
unit interval); we show that any two generalized arcs
are co-elementarily equivalent, and that co-elementary images of 
generalized arcs are generalized arcs.
We also show that if $f:X \to Y$ is a
function between compact Hausdorff spaces and if $X$ is an arc, then 
$f$ is a co-elementary map if and only if $Y$ is an arc and $f$ is a
monotone continuous surjection.  
\end{abstract}

\section{Introduction and Outline of Results.}\label{1}

A {\bf generalized arc\/} is a continuum (i.e., a connected compact Hausdorff
space) that has exactly two non-separating points; an {\bf arc\/} is a
metrizable generalized arc.  The class of generalized arcs is precisely the
class of linearly orderable continua; each generalized arc admitting exactly
two compatible linear orders.  The class of (continuous
images of) generalized arcs has been extensively 
studied over the years (see \cite{HY,NTT,Wil}); the most well-known results
in this area being that any two arcs are homeomorphic (to the standard closed
unit interval on the real line); and (Hahn-Mazurkiewicz) 
that a Hausdorff
space is a continuous image of an arc if and only if that space is a locally
connected metrizable continuum.  In this paper, a continuation of \cite{Ban3},
we study the model-theoretic topology of generalized arcs; in particular,
the ``dualized model theory'' of these spaces.  

Many notions from classical first-order model theory, principally elementary
equivalence and elementary embedding, may be phrased in terms of mapping
conditions involving the ultraproduct construction.  Because of the 
(Keisler-Shelah) ultrapower theorem (see, e.g., \cite{CK}), two relational
structures are elementarily equivalent if and only if some ultrapower of one
is isomorphic to some ultrapower of the other; a function from one relational
structure to another is an elementary embedding if and only if there is an
ultrapower isomorphism so that the obvious square mapping diagram commutes
(see, e.g., \cite{Ban2,Ban5,Ekl}).  The ultrapower construction in turn is a 
direct limit 
of direct products, and is hence capable of being transferred into a
purely category-theoretic setting.  In this paper we focus on the category
$\bf CH$
of compact Hausdorff spaces and continuous maps, but perform the transfer
into the opposite category (thus justifying the phrase
``dualized model theory'' above).

In $\bf CH$ one then constructs ultracoproducts,
and talks of co-elementary equivalence and co-elementary maps.  Co-elementary
equivalence is known \cite{Ban2,Ban5,Gur} to preserve important properties
of topological spaces, such as being infinite, being Boolean (i.e., totally
disconnected), having (Lebesgue) covering dimension $n$, and being a
decomposable continuum.  If $f:X \to Y$ is a co-elementary map in $\bf CH$,
then of course $X$ and $Y$ are co-elementarily equivalent (in symbols
$X \equiv Y$).  Moreover, since $f$ is a continuous surjection (see 
\cite{Ban2}), additional information about $X$ is transferred to $Y$.  For
instance, continuous surjections in $\bf CH$ cannot raise {\bf weight\/} (i.e., 
the smallest cardinality of a possible topological base, and for many
reasons the right cardinal invariant to replace cardiality in the dualized
model-theoretic setting), so metrizability
(i.e., being of countable weight in the compact Hausdorff context) is
preserved.  Also local connectedness is preserved, since continuous surjections
in $\bf CH$ are quotient maps.  Neither of these properties is an invariant
of co-elementary equivalence alone.  

When attention is restricted to the full subcategory of Boolean 
spaces, the dualized model theory matches perfectly with
the model theory of Boolean algebras because of Stone duality.  In the
larger category there is no such match \cite{Bana,Ros}, however, and one is
forced to look for other (less direct) model-theoretic aids.  Fortunately  
there is a finitely axiomatizable Horn class of bounded lattices,
the so-called {\it normal disjunctive\/} lattices \cite{Ban8} (also called
{\it Wallman\/} lattices in \cite{Ban5}), comprising precisely the (isomorphic
copies of) lattices that serve as bases for the closed sets of compact
Hausdorff spaces.  We go from lattices to spaces, as in the case of Stone
duality, via the {\bf maximal spectrum\/} $S(\;)$, pioneered by H. Wallman 
\cite{Walm}.
$S(A)$ is the space of 
maximal proper filters of $A$; a typical basic closed set in $S(A)$ is the
set of elements of $S(A)$ containing a given element of $A$.  $S(\;)$
is contravariantly functorial; if $f:A \to B$ is a homomorphism of normal
disjunctive lattices and $M \in S(B)$, then $f^S(M)$ is the unique maximal
filter in $A$ containing the pre-image of $M$ under $f$.  It is a fairly
straightforward task to show, then, that $S(\;)$ converts ultraproducts
to ultracoproducts, elementarily equivalent lattices to co-elementarily
equivalent compact Hausdorff spaces, and elementary embeddings to co-elementary
maps (see \cite{Ban2,Ban4,Ban5,Ban8,Gur}).  An important consequence of
this is a L\"{o}wenheim-Skolem theorem for co-elementary maps:  every 
compact Hausdorff space maps co-elementarily onto a compact metrizable space. 
(This result is used in \ref{2.4} and \ref{2.6} below.)

In \cite{Ban3} we showed that
any locally connected metrizable space co-elementarily equivalent to an arc
is already an arc; here we present the following results. 
$(i)$ if $f:X \to Y$ is a co-elementary
map in $\bf CH$, and if $Y$ is locally connected (in particular, a generalized
arc), then $f$ is a monotone continuous surjection; 
$(ii)$ co-elementary
images of (generalized) arcs are (generalized) arcs;
$(iii)$ any two generalized arcs are co-elementarily equivalent;
$(iv)$ if $X$ is a generalized arc and $f:X \to Y$ is
an irreducible co-elementary map in $\bf CH$, then $f$ is a homeomorphism; 
$(v)$ if every locally connected co-elementary pre-image of an arc is a
generalized arc, then every locally connected compact Hausdorff space
co-elementarily equivalent to a generalized arc is also a generalized arc;
and $(vi)$ if $X$ is an arc    
and $f$ is a function from $X$ to a compact Hausdorff space $Y$, then $f$ is
a co-elementary map if and only if $Y$ is an arc and $f$ is a monotone
continuous surjection.

Local connectedness is necessarily a part of $(v)$ above.
We do not know at present whether the hypothesis in $(v)$ is true; nor
do we know whether monotone surjections between 
generalized arcs are always co-elementary maps.  

\subsection{Remark.}\label{1.1}
By way of contrast, there is a Boolean analogue
to some of the results above.  Define a {\bf generalized Cantor set\/} to be any
non-empty Boolean space with no isolated points, and a {\bf Cantor set\/} to
be a metrizable generalized Cantor set.  It is well known that any two
Cantor sets are homeomorphic (to the standard Cantor middle thirds set in the
real line), and that the generalized Cantor sets are precisely the Stone
duals of the atomless Boolean algebras, constituting an elementary class 
whose first-order 
theory is $\aleph_0$-categorical, complete, and model complete.    
In $(ii)$ and $(iii)$, one may replace ``arc'' with ``Cantor set'' uniformly;
a straightforward application of $\aleph_0$-categoricity.	
The analog of $(iv)$ is
false (see Example 3.3.4$(iv)$ in \cite{Ban2}); the projective cover map to
a generalized Cantor set is always an irreducible co-elementary map between
(seldom-homeomorphic) generalized Cantor sets.  As for $(v)$, it follows from
the results on dimension in \cite{Ban2} that any compact Hausdorff space
co-elementarily equivalent to a generalized Cantor set is itself a 
generalized Cantor set.  Finally, regarding $(vi)$, {\it all\/} continuous 
surjections between
generalized Cantor sets are co-elementary maps.  This is a direct consequence
of the model completeness of the theory of atomless Boolean algebras.\\ 
  
\section{Methods and Proofs.}\label{2}
We begin with a proof of $(i)$ above.  Recall that a map $f:X \to Y$ is 
{\bf monotone\/} (resp. {\bf strongly monotone\/}) if the inverse image
of a point (resp. a closed connected subset) of $Y$ is connected in $X$.

\subsection{Proposition.}\label{2.1}
Let $f:X \to Y$ be a co-elementary
map in $\bf CH$, with $Y$ locally connected.  Then $f$ is a strongly
monotone continuous surjection.\\

\noindent
{\bf Proof.\/} Assume $f:X \to Y$ is co-elementary, $Y$ is locally
connected, and $f$ is not strongly monotone.  Then there is a subcontinuum 
$S$ of $Y$ such that
the inverse image $A := f^{-1}[S]$ is disconnected.  Since $A$ is closed,
we can write $A = A_1 \cup A_2$ where each $A_i$ is closed non-empty,
and $A_1 \cap A_2 = \emptyset$.  Let $U_i$ be an open neighborhood of $A_i$,
with $U_1 \cap U_2 = \emptyset$.  If $C$ is a subcontinuum  
of $X$ containing $A$, then we can pick 
some $x_C \in C \setminus (U_1 \cup U_2)$.  Let $B$ be the closure of
the set of all such points $x_C$, as $C$ ranges over all subcontinua
containing $A$.  Since no point $x_C$ lies in $U_1 \cup U_2$, $B$
is disjoint from $A$, but intersects every subcontinuum of $X$ that
contains $A$.  

Now $f[B]$ is closed in $Y$ and disjoint from $S$. Let $W$ be an open
neighborhood of $S$ whose closure misses $f[B]$.   
Since $Y$ is locally connected, we have, for each $y \in S$, a connected
open neighborhood $V_y$ of $y$ such that $V_y \subseteq W$.  Since $S$
is  connected, so also is $V := \bigcup_{y \in S}V_y$; and the closure
$K$ of $V$ is a subcontinuum containing $S$.  Since $V \subseteq W$, 
and the closure of $W$ is disjoint from $f[B]$, 
we know that $K$ is also disjoint from $f[B]$.  We need a fact proved 
elsewhere.\\

\noindent
{\bf Lemma.\/}(Lemma 2.8 in \cite{Ban5}) Let
$f:X \to Y$ be a co-elementary map in
$\bf CH$, with $K \subseteq Y$ a subcontinuum.  Then
there is a subcontinuum $C \subseteq X$ such that $K = f[C]$, and	
whenever $V \subseteq K$ is open in $Y$, $f^{-1}[V] \subseteq C$.\\

Using the Lemma, there exists a subcontinuum $C \subseteq X$ such that
$f[C] = K$ and $f^{-1}[V] \subseteq C$.  Let $x \in A$.  Then there
is a neighborhood $U$ of $x$ with $f[U] \subseteq V$.  Thus
$x \in U \subseteq f^{-1}[V] \subseteq C$, hence we infer $A \subseteq
C$.  Every subcontinuum of $X$ containing $A$ must intersect $B$, so  
$\emptyset \neq f[B \cap C] \subseteq f[B] \cap f[C] =
f[B] \cap K = \emptyset$.  This contradiction completes the proof. $\dashv$  

\subsection{Remark.}\label{2.2}
The Lemma above provides only a weak consequence of co-elementarity.  Indeed,
the usual projection map from the standard closed unit square in the plane
onto its first co\"{o}rdinate is not co-elementary because it does not
preserve topological dimension.  Nevertheless, it does satisfy the conclusion
of the Lemma.\\

Now we are in a position to prove $(ii)$.

\subsection{Proposition.}\label{2.3}
Let $f:X \to Y$ be a co-elementary map in $\bf CH$.  If $X$ is a generalized
arc, then so is $Y$.\\

\noindent
{\bf Proof.\/} Let $f:X \to Y$ be a co-elementary map in $\bf CH$, with
$X$ a generalized arc.  
$Y$ is a locally connected continuum
because $X$ is locally connected and $f$ is a continuous surjection.  By
\ref{2.1}, $f$ is monotone; it remains to show $Y$ has precisely two
non-separating points.

Let $a, b \in X$ be the two non-separating points of $X$.  $Y$ is 
non-degenerate because of co-elementarity; monotonicity then tells us that
$f(a) \neq f(b)$.  If $f(a)$ were to separate $Y$, we could also separate
$X \setminus K$, where $K:= f^{-1}[\{f(a)\}]$ is a subcontinuum 
(i.e., closed subinterval) containing
the endpoint $a$.  This is easily seen to be impossible for
generalized arcs.  Now let $y \in Y \setminus \{f(a),f(b)\}$, with 
$K := f^{-1}[\{f(y)\}]$.  Then $K$ is a subcontinuum of $X$ containing neither
endpoint.  Thus $X \setminus K$ is disconnected; hence $y$ separates $Y$.
We therefore conclude that $Y$ is a generalized arc. $\dashv$.\\

We can very quickly settle $(iii)$.

\subsection{Proposition.}\label{2.4}
Let $X$ and $Y$ be two generalized arcs.  Then $X \equiv Y$.\\

\noindent
{\bf Proof.\/} Let $X$ and $Y$ be generalized arcs.  By the 
L\"{o}wenheim-Skolem theorem for co-elementary maps, there exist co-elementary
maps $f:X \to X_0$ and $g:Y \to Y_0$, where $X_0$ and $Y_0$ are compact
metrizable.  By \ref{2.3}, the images are generalized arcs; hence they are
arcs.  Thus $X_0$ and $Y_0$ are homeomorphic, and we conclude $X \equiv Y$
because \cite{Ban2} co-elementary equivalence is an honest equivalence
relation. $\dashv$\\

To handle $(iv)$, recall that a continuous surjection $f:X \to Y$ is
{\bf irreducible\/} if $Y$ is not the image under $f$ of a proper closed
subset of $X$.

\subsection{Proposition.}\label{2.5}
Let $f:X \to Y$ be an irreducible co-elementary map in $\bf CH$.  If $X$
is a generalized arc, then $f$ is a homeomorphism.\\

\noindent
{\bf Proof.\/}  It suffices to show $f$ is one-one.  Let $y \in Y$, with
$K := f^{-1}[\{y\}]$, a subcontinuum of $X$ by \ref{2.1}.  Since $X$ is
a generalized arc, $K$ is either a singleton or a closed subinterval with
non-empty interior.  The latter case easily contradicts the irreducibility
of $f$, however. $\dashv$\\

In \cite{Gur} it is shown that every infinite compact Hausdorff space is 
co-elementarily equivalent to a compact Hausdorff space that is not locally
connected.  (See also \cite{Ban5} for refinements.)  This explains the
necessity of the local connectedness hypothesis in $(v)$.  

\subsection{Proposition.}\label{2.6}      
Suppose every locally connected co-elementary pre-image of an arc is a
generalized arc.  Then every locally connected compact Hausdorff space
co-elementarily equivalent to a generalized arc is itself a generalized arc.\\

\noindent
{\bf Proof.\/}  Suppose $X \in \bf CH$ is locally connected,
$X \equiv Y$, and $Y$ is a generalized arc.  As in the proof of \ref{2.4}
above, we have co-elementary maps $f:X \to X_0$ and $g:Y \to Y_0$, where
$X_0$ and $Y_0$ are metrizable.  Furthermore, we know that $X_0$ is 
locally connected and that $Y_0$ is an arc (\ref{2.1} again).  By the
transitivity of co-elementary equivalence, we know $X_0 \equiv Y_0$; by
the main result of \cite{Ban3}, we know $X_0$ is an arc.  Our hypothesis
then tells us that $X$ is a generalized arc. $\dashv$\\

We finish with a proof of $(vi)$.  If $X$ is an arc and $f:X \to Y$ is 
a co-elementary map in $\bf CH$, then $Y$ is an arc and $f$ is a monotone
continuous surjection by \ref{2.1} and \ref{2.2}.  So it suffices to prove
the following. 

\subsection{Proposition.}\label{2.7}
Every monotone continuous surjection from
an arc to itself is a co-elementary map.\\

\noindent
{\bf Proof. } Let us take our arc to be the standard closed unit interval
${\bf I}$ with its usual order.
$f$ is either $\leq$-preserving or $\leq$-reversing, so we lose
no generality in assuming $f$ to be the former.

For any topological space $X$, we denote the closed set lattice of $X$ by
$F(X)$.  $F(\;)$ converts continuous maps contravariantly into lattice
homomorphisms, and serves as a right inverse for $S(\;)$: $S(F(X))$ is 
naturally homeomorphic to $X$ for any compact Hausdorff $X$.
  
Monotone continuous surjections from $\bf I$ to itself are strongly monotone; 
hence
$f^F:F({\bf I}) \to F({\bf I})$ is a lattice embedding that takes closed
intervals (in this case the connected elements of the lattice) to closed 
intervals.  
However, $f^F$ will take atoms to non-atoms when $f$ is not injective.
Thus $f^F$ is not an elementary embedding without being an isomorphism.
The idea is to restrict the domain and range of $f^F$ in such a way that
the resulting lattice embedding, call it $g$, is elementary, and $g^S = f$.

Our plan is to create an elementary lattice embedding $g:{\cal A} \to
{\cal B}$, where ${\cal A}$ and ${\cal B}$ are atomless lattice bases
for ${\bf I}$ (i.e., both $\cal A$ and $\cal B$ are atomless, as well as 
meet-dense in $F({\bf I})$), and $g$ agrees with the restriction 
of $f^F$ to ${\cal A}$.

Since $S(\cal A)$ and $S(\cal B)$ are naturally homeomorphic to $\bf I$, 
and $f$ is just $g^S$ conjugated with these homeomorphisms, $f$ is a
co-elementary map provided $g^S$ is.

For each $y \in \bf I$, let $\lambda (y) := \mbox{inf}(f^{-1}[\{y\}])$
and $\rho (y) := \mbox{sup}(f^{-1}[\{y\}])$.  Then for any closed interval
$[x,y] \in F(\bf I)$, $f^F([x,y]) = [\lambda (x),\rho (y)]$.  Both $\lambda$
and $\rho$ are right inverses for $f$, and are hence strictly increasing
(but not necessarily continuous).  Of course $\lambda (0)  = 0$ and
$\rho (1) = 1$.  

Let $L,R \subseteq \bf I$, with $0 \in L$ and $1 \in R$.  If ${\cal I}(L,R)$
denotes the set of all finite unions of intervals $[x,y]$ with $x \in L$
and $y \in R$, then ${\cal I}(L,R)$ is a sublattice of $F(\bf I)$, which is
atomless just in case $L \cap R = \emptyset$.  If $L$ and $R$ are dense
in $\bf I$, then ${\cal I}(L,R)$ is a lattice base as well.

Now fix $L,R \subseteq \bf I$ to be disjoint countable dense subsets, 
with $0 \in L$ and $1 \in R$, and set ${\cal A} := {\cal I}(L,R)$.
Then the image of $\cal A$ under $f^F$  is ${\cal I}(\lambda [L],\rho[R])$.
Clearly $\lambda [L] \cap \rho [R] = \emptyset$, $0 \in \lambda [L]$, and
$1 \in \rho [R]$.  
Let $L', R' \subseteq \bf I$ be disjoint countable dense
subsets, with $\lambda [L] \subseteq L'$, $\rho [R] \subseteq R'$, and
set ${\cal B} := {\cal I}(L',R')$.   Then $\cal B$ is a countable atomless
lattice base for $F(\bf I)$, and we denote by $g:{\cal A} \to \cal B$
the embedding $f^F$ with its domain and range so restricted.  It remains
to show that $g$ is an elementary embedding, and for this it suffices to
show that for
each finite set $S$ in $\cal A$ and each $b \in \cal B$, there is an
automorphism on $\cal B$ that fixes $g[S]$ pointwise
and takes $b$ into $g[\cal A]$. 

Let $x_1, ..., x_n$ be a listing, in increasing order, of the endpoints
of the component intervals of $g[S] \cup \{b\}$ (so each $x_i$ is in
$L' \cup R'$), with $X_i := 
f^{-1}[\{f(x_i)\}]$, $1 \leq i \leq n$.  Each $X_i$ is either a singleton
or a non-degenerate closed interval, and 
for $1 \leq i < j \leq n$,
either $X_i = X_j$ or each element of $X_i$ is less than each element of
$X_j$.  Let $U_i$ be an open-interval neighborhood of $X_i$ such that
$U_i \cap U_j = \emptyset$ whenever $X_i \neq X_j$.  Since $f$ is a
$\leq$-preserving surjection and the sets $L$ and $R$ are dense in $\bf I$,
each $U_i$ has
infinite intersection with both $\lambda [L]$ and $\rho [R]$.  If
$x_i \in \lambda [L] \cup \rho [R]$, set $x_i' := x_i$.  Otherwise we
know $x_i$ is an endpoint of a component interval of $b$; and we choose
$x_i' \in U_i$ in such a way that $x_i' \in \lambda [L]$ if and only if
$x_i \in L'$, and $x_i' < x_j'$ whenever $x_i < x_j$ and
$X_i = X_j$.  This 
procedure produces an increasing sequence $x_1',...,x_n'$ of elements
of $\lambda [L] \cup \rho [R]$; $x_i' \in \lambda [L]$ if and only if
$x_i \in L'$.  For each $a \in g[S] \cup \{b\}$, let $a'$ be built up using 
the endpoints $x_i'$ in the same way as $a$ is built up using the endpoints
$x_i$.  Then $a' = a$ for each $a \in g[S]$, and $b' \in g[\cal A]$.
Now by a classic (Cantor) back and forth argument, there is an order
automorphism on $L' \cup R'$ that fixes $L'$ and $R'$ setwise and takes
$x_i$ to $x_i'$ for $1 \leq i \leq n$.  This order automorphism gives
rise to the lattice automorphism on $\cal B$ that we require. $\dashv$\\

\end{document}